\documentclass[12pt]{article}
\usepackage[dvips]{epsfig}
\usepackage{amsthm,amsmath,amssymb}
\usepackage{latexsym}
\usepackage{graphics}
\usepackage{url}
\usepackage{placeins}
\usepackage{hyperref}
\usepackage{color}
\usepackage{bm}
\usepackage{comment}
\usepackage{statex}

\definecolor{blue}{rgb}{0,0,0.9}
\definecolor{red}{rgb}{0.9,0,0}
\definecolor{green}{rgb}{0,0.50,0.10}
\definecolor{violet}{rgb}{0.5804,0.0000,0.8275}

\def\@themcountersep{}

\newtheorem{THEO}{Theorem}[section]
\newtheorem{ALGo}[THEO]{Algorithm}
\newtheorem{CONJ}[THEO]{Conjecture}
\newtheorem{COND}[THEO]{Condition}
\newtheorem{CORO}[THEO]{Corollary}
\newtheorem{DEFI}[THEO]{Definition}
\newtheorem{EXAMP}[THEO]{Example}
\newtheorem{FACT}[THEO]{Fact}
\newtheorem{HYPO}[THEO]{Hypothesis}
\newtheorem{LEMM}[THEO]{Lemma}
\newtheorem{PROB}[THEO]{Problem}
\newtheorem{PROP}[THEO]{Proposition}
\newtheorem{REMA}[THEO]{Remark}
\newcommand{\theo}{\begin{THEO}}
\newcommand{\algo}{\begin{ALGo} \rm}
\newcommand{\cond}{\begin{COND}}
\newcommand{\conj}{\begin{CONJ}}
\newcommand{\coro}{\begin{CORO}}
\newcommand{\defi}{\begin{DEFI} \rm}
\newcommand{\examp}{\begin{EXAMP} \rm}
\newcommand{\fact}{\begin{FACT}}
\newcommand{\hypo}{\begin{HYPO} \rm}
\newcommand{\lemm}{\begin{LEMM}}
\newcommand{\prob}{\begin{PROB} \rm}
\newcommand{\prop}{\begin{PROP}}
\newcommand{\rema}{\begin{REMA} \rm}
\newcommand{\etheo}{\end{THEO}}
\newcommand{\ealgo}{\end{ALGo}}
\newcommand{\econd}{\end{COND}}
\newcommand{\econj}{\end{CONJ}}
\newcommand{\ecoro}{\end{CORO}}
\newcommand{\edefi}{\end{DEFI}}
\newcommand{\eexamp}{\end{EXAMP}}
\newcommand{\efact}{\end{FACT}}
\newcommand{\ehypo}{\end{HYPO}}
\newcommand{\elemm}{\end{LEMM}}
\newcommand{\eprob}{\end{PROB}}
\newcommand{\eprop}{\end{PROP}}
\newcommand{\erema}{\end{REMA}}

\def\0{\mbox{\bf 0}}
\def\1{\mbox{\bf 1}}
\def\2{\mbox{\bf 2}}
\def\3{\mbox{\bf 3}}
\def\4{\mbox{\bf 4}}
\def\5{\mbox{\bf 5}}
\def\6{\mbox{\bf 6}}
\def\7{\mbox{\bf 7}}
\def\8{\mbox{\bf 8}}
\def\9{\mbox{\bf 9}}

\def\e{\mbox{\boldmath $e$}}

\def\g{\mbox{\boldmath $g$}}

\def\s{\mbox{\boldmath $s$}}

\def\u{\mbox{\boldmath $u$}}
\def\v{\mbox{\boldmath $v$}}

\def\x{\mbox{\boldmath $x$}}

\def\G{\mbox{\boldmath $G$}}

\def\Q{\mbox{\boldmath $Q$}}
\def\R{\mbox{\boldmath $R$}}

\def\U{\mbox{\boldmath $U$}}
\def\V{\mbox{\boldmath $V$}}

\def\X{\mbox{\boldmath $X$}}

\def\IC{\mbox{$\cal I$}}
\def\JC{\mbox{$\cal J$}}

\def\inprod#1#2{\langle#1, \, #2\rangle}

\def\Real{\mbox{$\mathbb{R}$}}
\def\Integer{\mbox{$\mathbb{Z}$}}

\def\s0{\mbox{\scriptsize \boldmath $0$}}

\def\bPhi{\mbox{\boldmath $\Phi$}}

\def\Real{\mathbb{R}}
\def\coneK{\mathbb{K}}

\def\spaceV{\mathbb{V}}
\def\SymMat{\mathbb{S}}

\def\Integer{\mathbb{Z}}

\begin{document}


\title{ \Large 
Generating Cutting Inequalities Successively
\\  for Quadratic Optimization Problems 
in Binary Variables \\
}

\author{
\normalsize 
Sunyoung Kim\thanks{Department of Mathematics, Ewha W. University, Seoul, 52 Ewhayeodae-gil, Sudaemoon-gu, Seoul 120-750, Korea 
			({\tt skim@ewha.ac.kr}). 
The research was supported  by   NRF 2021-R1A2C1003810.
}, \and \normalsize
Masakazu Kojima\thanks{Department of Industrial and Systems Engineering,
	Chuo University,  Tokyo 192-0393, Japan 
	 ({\tt kojima@is.titech.ac.jp}).
	 This research was supported by Grant-in-Aid for Scientific Research (A) 19H00808.
}
}


\date{\normalsize\today}

\maketitle


\begin{abstract}
\noindent
We propose a successive generation of cutting inequalities for binary quadratic optimization problems.
Multiple cutting inequalities are successively generated for the convex hull of  the set of the  optimal solutions $\subset \{0, 1\}^n$,
while the standard  cutting inequalities are used for the convex hull of the feasible region. 
An arbitrary linear inequality with integer coefficients and the right-hand side value in integer is considered as a candidate
for a valid inequality.  The validity of the linear inequality 
is determined by solving a conic relaxation of a subproblem such as the doubly nonnegative relaxation,
under the assumption that an upper bound  for the unknown optimal value of the problem is available.
Moreover,  the candidates generated for the multiple cutting inequalities 
are tested simultaneously for their validity in parallel.
Preliminary numerical results on 60 quadratic unconstrained binary optimization problems with a simple implementation of the successive
 cutting inequalities using an 8- or 32-core machine show that 
 the exact optimal  values are obtained for 70\%
of the tested problems, demonstrating the strong potential of the proposed technique.
 \end{abstract}

\vspace{0.5cm}


\vspace{0.5cm}



\section{Introduction}

Cutting inequalities \cite{BALAS1993,BONAMI2019,ENGAU2010,HELMBERG1998,MARCHAND2002} and
conic relaxations such as linear programming (LP), semidefinite programming (SDP) \cite{ANJOS2012}, 
doubly nonnegative (DNN) \cite{KIM2013,ARIMA2017} have been regarded as the two most basic tools
for solving nonconvex and/or combinatorial optimization problems.
They have been frequently 
incorporated in the brach-and-bound and branch-and-cut framework \cite{GLOVER1998,GUIMARAES2020,KRISLOCK2014,KRISLOCK2017}
to solve the problems.

The main purpose of this paper is to propose 
a {\em successive cutting inequality technique}, abbreviated by {\em SCIT},  for binary quadratic optimization problems (QOPs, {\it i.e.}, QOPs in binary variables), and
to demonstrate its strong potential to become a very powerful tool for solving binary QOPs, 
 through preliminary numerical results 
 by an experimental method that implements the very basics of SCIT.

To describe the motivation and  basic idea of SCIT, we consider a 
general nonconvex optimization problem:  
\[ 
\mbox{P: } \zeta = \min \{f(x) : \x \in S \}, 
\] 
where $f$ denotes a real valued function on the $n$-dimensional Euclidean space $\Real^n$ and 
$S$ denotes a closed subset of $\Real^n$. 
The assumption that $f$ is a polynomial function in $\x \in \Real^n$ and $S$ is described 
by polynomial equalities and inequalities \cite{LASSERRE2001}  is required at least  
 for the discussion of 
a conic relaxation of problem P.
While the case where those polynomials are linear or quadratic is mainly dealt with,
such assumptions are not so relevant in the discussion below. 

Except for LP relaxation, 
a conic relaxation problem with a linear objective function over a 
closed convex feasible region $\widehat{S}$ is embedded in a different space $\spaceV$, 
often called a lifted space,  with a higher dimension such as 
the linear space of symmetric matrices. $\widehat{S}$ is described by linear equalities and inequalities 
in $\spaceV$ and a closed convex cone $\coneK \subset \spaceV$. 
In short, 
when the cone $\coneK$ used is the nonnegative orthant of the Euclidean space, the positive semidefinite matrix cone or the doubly nonnegative 
matrix cone, the conic relaxation is called an LP relaxation, an SDP relaxation or a DNN relaxation, 
respectively.   
The lifted space $\spaceV$ is identified 
with the original $n$-dimensional Euclidean space $\Real^n$ where problem P is defined,  
and $f$ itself is assumed to be linear, for simplicity of discussion below.
Then, it is clear that 
$\zeta = \min \{ f(\x) : \mbox{co}(S) \}$, where $\mbox{co}( S)$ denotes the convex hull of 
$S$. 

Under the above setting and assumptions, the conic relaxation of~P can be described as 
a convex optimization problem $\widehat{\rm P}$ with a linear objective function 
: $\min \hat{\zeta} = \{ \hat{f}(\x) : \x \in \widehat{S} \}$, where $\widehat{S}$ is a 
closed convex subset of $\spaceV$ (identified with $\Real^n$) such that 
$S \subseteq \mbox{co}(S) \subseteq \widehat{S}$ 
and $\hat{f}(\x) = f(\x)$ for every $\x \in S$. 
In practice, $\widehat{\rm P}$ is constructed as a numerically tractable convex optimization problem. 
If $\mbox{co}(S) = \widehat{S}$ held, we would have the optimal value $\zeta$ of problem~P 
by solving the relaxation problem $\widehat{\rm P}$. 
In most cases, however, 
 $\mbox{co}(S)$ is a proper subset 
of $\widehat{S}$, $\hat{\zeta} < \zeta$, and an optimal solution $\bar{\x}$ of  $\widehat{\rm P}$ 
is not a feasible solution of P. 
As a result, 
 an inequality which cuts off $\bar{\x}$ from $\widehat{S}$ 
but does not remove any $\x$ from $\mbox{co}(S)$ is desired to improve the lower bound 
$\hat{\zeta}$ of $\zeta$ and to 
compute a feasible approximate optimal solution of~P. This is a standard role and usage of cutting inequalities. 
Note that a cutting inequality is chosen from the family of valid inequalities of $S$.
In fact, $\mbox{co}(S)$ 
can be described as the set of points $\x$ which satisfies all the valid inequalities of $S$. 
Well-known 
triangle inequalities, which forms a sub-family of the valid inequalities 
for the binary polytopes, are frequently used to strengthen the SDP relaxation of binary QOPs  \cite{KRISLOCK2014}. 

In this paper,  
a cutting inequality plays a more active role 
under the additional assumption that $S \subseteq \{0,1\}^n$. 
We propose to generate 
a cutting inequality for the convex hull $\mbox{co}(S^*)$ of the set $S^*$ of optimal solutions 
of P, an inequality which is aimed at cutting off $\x \not\in \mbox{co}(S^*)$ from $\mbox{co}(S^*)$ 
(Recall that the standard cutting inequality is for the convex hull of the feasible region $S$). 
Generating such a cutting inequality is based on the following ideas:  Assume that an upper bound 
$\eta$ of the unknown optimal value $\zeta$ is available. Let $\g$ be an arbitrary integer column vector 
 in $\Real^n$. For every integer $\alpha$, we consider a pair of subset $S(\alpha) = \{\x \in S : \g^T \x \leq \alpha\}$ 
and $S(\alpha)^+ = \{\x \in S : \g^T \x \geq \alpha + 1\}$. Since $\g^T\x^*$ is an integer  for every 
$\x^* \in S^* \subset \{0,1\}^n$, it is obvious that $S^*$ is included in the union of $S(\alpha)$ 
and $S(\alpha)^+$. 
Hence if $S^* \cap S(\alpha) = \emptyset$, then $S^* \subseteq S(\alpha)^+$, {\it i.e.}, 
$\g^T \x \geq \alpha + 1$ 
serves as a cutting inequality for $S^*$. To obtain a certificate of  $S^* \cap S(\alpha) = \emptyset$, 
we solve a conic relaxation  $\widehat{\rm P}(\alpha)$: 
$\hat{\zeta}(\alpha)= \min\{\hat{f}(\x) : \x \in \widehat{S}(\alpha)\}$ of a subproblem ${\rm P}(\alpha)$: 
$\zeta(\alpha) = \min\{f(\x) : \x \in S(\alpha)\}$, 
where $\widehat{S}(\alpha)$ denotes a closed convex subset of 
$\spaceV = \Real^n$ containing $S(\alpha)$. If $\eta < \hat{\zeta}(\alpha)$ holds, 
then $\zeta \leq \eta < \hat{\zeta}(\alpha)$; 
hence $S(\alpha)\cap S^* \subseteq \widehat{S}(\alpha)\cap S^* = \emptyset$. Therefore, 
the inequality $\eta < \hat{\zeta}(\alpha)$ is a certificate for 
 $\g^T \x \geq \alpha + 1$ to be a 
cutting inequality for $S^*$. 
The largest $\alpha$ such that $\eta < \hat{\zeta}(\alpha)$ is most desirable 
to cut off a larger portion of $S \backslash S^*$. 
For such an $\alpha$, a $1$-dimensional search can be applied 
over the set of integers with starting $\alpha = -1$ since $\g^T \x \geq 0$ is a trivial cutting inequality for $S^*$. 

In the proposed SCIT, 
multiple candidates for cutting inequalities in the form $\g_j^T\x \geq \alpha_j + 1$ 
$(j=1,\ldots,m)$ are arranged before the iteration starts and set $\alpha_j + 1 = 0$ so that 
$\g_j^T\x \geq 0$  becomes a trivial cutting inequality for $S^* \ (j=1,\ldots,m)$. 
At each iteration of SCIT,  it verifies whether 
$\g_j^T\x \geq \alpha'_j + 1$ remains a valid inequality for some $\alpha_j' > \alpha_j$ 
by solving a conic relaxation problem for all $j=1,\ldots,m$. If it does, 
then $\alpha_j$ is updated to $\alpha_j'$, 
otherwise $\alpha'_j$ is replaced by a smaller 
$\alpha'_j \in [\alpha_j,\alpha_j')$ for the next iteration.  
Notice that these verifications and updates with $j=1,\ldots,m$ can be simultaneously performed  (within one iteration)  
in parallel. 

The effectiveness and efficiency of the performance of  SCIT
on large scale binary QOPs 
is dependent on 
the followings: \vspace{-1.6mm} 
\begin{description}
\item{(I) } A tight upper bound $\eta$ is available for the unknown optimal value $\zeta$. \vspace{-2mm} 
\item{(II) } A strong conic relaxation method that generates a tight lower bound $\hat{\zeta}(\alpha)$ for 
$\zeta(\alpha)$ can be utilized. \vspace{-2mm} 
\item{(III) } A powerful computer system can be used for parallel computing. \vspace{-2mm} 
\end{description}

As an application of SCIT,
we consider quadratic unconstrained binary 
optimization problems (QUBOs)  in Section 4. 
There exist many heuristic methods, which can be used for computing a tight 
upper bound $\eta$ of the optimal value $\zeta$ of a QUBO,
such as the tabu search \cite{GLOVER1998} and 
the genetic algorithm \cite{MITCHELL1998}. 
For (II), we utilize the Lagrangian-DNN relaxation \cite{KIM2013,ARIMA2017,ARIMA2018}, 
which is known to be much stronger than the standard SDP relaxation, and 
NewtBracket \cite{KIM2020} (the Newton-bracketing method \cite{KIM2019a}) as a numerical method to compute its optimal value. For (III), preliminary 
numerical results with a small scale computer are reported. 
 Inequalities 
of the form $\sum_{i \in I} x_i \geq \alpha + 1$ for some $I \subseteq \{1,\ldots,n\}$ are considered as the candidates for the cutting inequalities.
 If the inequality is shown to be 
a cutting inequality with $I=\{j\}$ and $\alpha +1  = 1$, then 
$x_j$ can be fixed to $x_j = 1$  and the size of the QUBO to be 
solved can be reduced. This is an important feature of SCIT. 

We investigate the numerical performance of SCIT  through an experimental method   on  
 $60$ QUBO instances with dimensions up to $250$ from BIQMAC \cite{BIQMAC}. 
Although the method is just a 
simple implementation of SCIT, not a well-designed software for solving QUBOs, 
it attained the exact optimal value within $10$ iterations 
for $70$\% cases of the $60$ instances. This is a remarkable result, which could not be
 expected.
It shows the promising potential of SCIT when it is incorporated into the branch-and-bound  method \cite{GLOVER1998,GUIMARAES2020,KRISLOCK2014,KRISLOCK2017}. 
We mention that theoretical aspects of SCIT including the convergence to 
the convex hull of $S^*$ are not dealt with here.

In Section 2, we present the fundamental facts which our construction of 
cutting inequalities build on after introducing notation and 
symbols. We present some details on SCIT in Section~3, 
and discuss its application to QUBOs in Section~4. 
The preliminary numerical results mentioned above are given in Section 4.4.
We conclude in Section 5.

\section{Preliminaries}

\subsection{Notation and symbols}

Let 
$ 
\Real = \mbox{the set of real numbers} \ \mbox{ and } 
\Integer = \mbox{the set of integers}.
$ 
For $\spaceV = \Real$  or 
$\Integer$, 
$\spaceV^n$ denotes the set of $n$-dimensional column vectors 
$(v_1,\ldots,v_n)$ with elements $v_i \in \spaceV$ $(i=1,\ldots,n)$, and 
$\spaceV^{\ell \times \ell}$ the set of $\ell \times \ell$ 
matrices $\V = [V_{ij}]$ with elements $V_{ij} \in \spaceV$ $(1\leq i,j \leq \ell)$. 
In particular, $\Real^n$ denotes the $n$-dimensional 
Euclidean space of column vectors. 
$\v^T$ stands for  
the transposed row vector of $\v$ for every $\v \in \spaceV^n$, and  $\u^T\v$ the 
inner product $\sum_{i=1}^n u_iv_i$ of $\u, \ \v \in \spaceV^n$. 
For $\U, \V \in \spaceV^{\ell \times \ell}$, their inner product is written as
$\inprod{\U}{\V} = \sum_{i=1}^{\ell}\sum_{j=1}^{\ell} U_{ij}V_{ij}$. 
Let 
\begin{eqnarray*}
\SymMat^{\ell} & = & \mbox{the linear space of $\ell \times \ell$ symmetric 
matrices $\X = [X_{ij}]$ $(1 \leq i,j \leq \ell)$}, \\
\SymMat^{\ell}_+ & = & \mbox{the cone of positive semidefinite matrices in $\SymMat^{\ell}$}. 
\end{eqnarray*}
Throughout the paper, $\widehat{\ }$ and $\hat{\ }$ (also $\widetilde{\ }$ and $\tilde{\ }$) are used for conic relaxation problems such that 
a conic relaxation problem $\widehat{\rm P}$~\eqref{eq:ConicRelaxation0} of the optimization 
problem P~\eqref{eq:NLP0} and their optimal values $\hat{\zeta}$ and $\zeta$, respectively. We use the subscripts $_{\rm a}$ and 
$_{\rm b}$ for a pair of subproblems 
obtained from 
their common parent problem by adding a cutting inequality; 
for example, a pair of subproblems P$_{\rm a}$~\eqref{eq:NLPa} and P$_{\rm a}$~\eqref{eq:NLPb} of 
P~\eqref{eq:NLP0} and their optimal values $\zeta_{\rm a}$ and $\zeta_{\rm b}$, respectively. 

\subsection{Basic ideas to generate cutting inequalities}

We begin with the following simple facts on which our cutting 
inequalities are constructed. 

\lemm \label{lemma1}
Let $\zeta, \ \zeta_{\rm a}, \ \zeta_{\rm b}, \ \eta, \
 \hat{\zeta}_{\rm a}$, and $\hat{\zeta}_{\rm b}$ be  real numbers satisfying 
$ 
\zeta = \min\{\zeta_{\rm a},\zeta_{\rm b} \} \leq  \eta, 
\ \hat{\zeta}_{\rm a} \leq \zeta_{\rm a} \ \mbox{ and } \hat{\zeta}_{\rm b} 
\leq \zeta_{\rm b}.
$ 
Assume that $\eta < \hat{\zeta}_{\rm a}$. Then $\hat{\zeta}_{\rm b} \leq \zeta_{\rm b} = \zeta$. 
\elemm
\proof{It follows from $\zeta = \min \{ \zeta_{\rm a}, \zeta_{\rm b}\}$ that 
at least one of $\zeta=\zeta_{\rm a}$ and $\zeta=\zeta_{\rm b}$ holds. 
If $\zeta_{\rm a} = \zeta$ held, then we would have $\hat{\zeta}_{\rm a} 
\leq \zeta_{\rm a} = \zeta \leq \eta$. This contradicts to the assumption that 
$\eta < \hat{\zeta}_{\rm a}$. 
\qed}

\bigskip 

\noindent
We note that the cutting inequalities  in Section 3.3 are constructed by Lemma \ref{lemma1}. More precisely, 
 \vspace{-2mm}
\begin{itemize}
\item $\zeta$ corresponds to the unknown optimal (minimum) value of the optimization 
problem ${\rm P}$~\eqref{eq:NLP0}, our target problem to solve, 
and $\eta$ to a known upper bound of 
$\zeta$.   \vspace{-2mm}
\item $\zeta_{\rm a}$ and $\zeta_{\rm b}$ correspond to optimal values of a pair of 
subproblems P$_{\rm a}$~\eqref{eq:NLPa} and P$_{\rm b}$~\eqref{eq:NLPb}, 
which are generated by adding cut inequalities to the feasible 
region of ${\rm P}$~\eqref{eq:NLP0}; hence $\zeta \leq \zeta_{\rm a}$ and 
$\zeta \leq \zeta_{\rm b}$. The identity 
$\min\{\zeta_{\rm a},\zeta_{\rm b} \} =  \zeta$ means at least one 
of P$_{\rm a}$ and P$_{\rm b}$ attains the same objective value as ${\rm P}$~\eqref{eq:NLP0}. 
\item $\hat{\zeta}_{\rm a}$ and $\hat{\zeta}_{\rm b}$ correspond to the optimal values of  
$\widehat{\rm P}_a$~\eqref{eq:ConicRelaxationa}  and $\widehat{\rm P}_{\rm b}$~\eqref{eq:ConicRelaxationb}, which are 
conic relaxations of P$_{\rm a}$~\eqref{eq:NLPa} and P$_{\rm b}$~\eqref{eq:NLPb}, 
respectively. 
\end{itemize}

\lemm \label{lemma2}
Let $\zeta, \ \eta, \ \hat{\zeta}, \ \tilde{\zeta}_{\rm a}$ and $\tilde{\zeta}_{\rm b}$ be  real numbers satisfying 
$ 
\hat{\zeta} \leq  \min\{\tilde{\zeta}_{\rm a} \ \mbox{and } \tilde{\zeta}_{\rm b} \} \leq \zeta \leq \eta.
$ 
Assume that $\eta < \tilde{\zeta}_{\rm a}$. Then, $\hat{\zeta} \leq \tilde{\zeta}_{\rm b} 
\leq \zeta$. 
\elemm
\proof{ Obvious. }

\bigskip 

\noindent
By Lemma \ref{lemma2}, 
the cutting inequalities in Section 3.4 are constructed, where  
\begin{itemize}
\item $\zeta$ corresponds to the unknown optimal (minimum) value of the optimization 
problem P~\eqref{eq:NLP0}, the target problem, and $\eta$ to a known upper bound of 
$\zeta$. \vspace{-2mm}
\item $\hat{\zeta}$ corresponds to the 
optimal values of a conic relaxation problem $\hat{\rm P}$~\eqref{eq:ConicRelaxation0} 
of P~\eqref{eq:NLP0}. \vspace{-2mm}
\item $\tilde{\zeta}_{\rm a}$ and $\tilde{\zeta}_{\rm b}$ correspond the optimal values of a pair of
 subproblems $\widetilde{\rm P}_{\rm a}$~\eqref{eq:ConicRelaxationA} and 
$\widetilde{\rm P}_{\rm b}$~\eqref{eq:ConicRelaxationB} of $\widehat{\rm P}$~\eqref{eq:ConicRelaxation0}, 
which are generated by adding cutting inequalities to 
the feasible region of $\widehat{\rm P}$~\eqref{eq:ConicRelaxation0}. The inequality 
$\min\{\tilde{\zeta}_{\rm a},\tilde{\zeta}_{\rm b}\}  \leq \zeta$ means at least 
one of $\widetilde{\rm P}_{\rm a}$~\eqref{eq:ConicRelaxationA} and 
$\widetilde{\rm P}_{\rm b}$~\eqref{eq:ConicRelaxationB}  
acts as a conic relaxation 
of P~\eqref{eq:NLP0}. 
\end{itemize}

\section{Conic relaxations of optimization problems in binary variable with cutting inequalities}

\subsection{An optimization problem in binary variables} 

Throughout this section, we consider the following nonconvex optimization problem in 
binary variables $x_i \in \{0,1\}$ $(i=1,\ldots,n)$: 
\begin{eqnarray}
\mbox{P: } \zeta & = & \min \left\{f(\x) : \x \in S \right\}, \label{eq:NLP0} 
\end{eqnarray}
where
\begin{description}
\item{(a) } $\emptyset \not= S \subset \{0,1\}^{n}$,
\item{(b) } $f(\x) \in \Integer$ for every $\x \in \{0,1\}^{n}$. 
\end{description} 
Under these conditions, problem P has an optimal solution $\x^*$. 
The optimal value $\zeta$ and solution $\x^*$ are unknown. Moreover, 
we assume  that  
\begin{description}
\item{(c) } An upper bound $\eta \in \Integer$  for $\zeta$ is available. 
\end{description} 

\subsection{A conic relaxation of problem P in the lifted symmetric matrix space}

We first introduce a finite dimensional vector space 
into which 
the conic relaxation is embedded. For simplicity of discussion and convenience 
of presenting an application in Section~4, 
we focus on the case where the linear space is $\SymMat^{\ell}$ of 
$\ell \times \ell$ symmetric matrices.

Consider 
\begin{eqnarray}
\mbox{$\widehat{\rm P}$: } \hat{\zeta} & = & \min \left\{ \inprod{\Q}{\X} : \X \in \widehat{S} \right\}. \label{eq:ConicRelaxation0}
\end{eqnarray}
Here \vspace{-1mm} 
\begin{description}
\item{(d) } $\Q$ is a matrix in $\SymMat^{\ell} \cap \Integer^{\ell\times \ell}$ such that 
$\inprod{\Q}{\bPhi(\x)} = f(\x)$ for every $\x \in \{0,1\}^n$. \vspace{-2mm} 
\item{(e) } $\bPhi$ is a mapping from $\Real^n$ into $\SymMat^{\ell}$ such that 
$\bPhi(\x) \in \SymMat^{\ell} \cap \Integer^{\ell \times \ell}$ if 
$\x \in \{0,1\}^n$.\vspace{-2mm} 
\item{(f) } $\widehat{S}$ is a closed convex subset of $\SymMat^{\ell}$ such that 
$\bPhi(S) \subset \widehat{S}$, {\it i.e.}, 
$\bPhi(\x) \in \widehat{S}$ for every $\x \in S$.\vspace{-1mm} 
\end{description}
These three conditions characterize problem $\widehat{\rm P}$~\eqref{eq:ConicRelaxation0} as a conic  
(SDP and DNN) relaxation problem of P~\eqref{eq:NLP0} in the space $\SymMat^{\ell}$. 
In particular, 
$ 
\hat{\zeta} \leq \zeta \leq \eta \ \mbox{and } \bPhi(\x^*) \in \widehat{S}. 
$ 

\examp \label{example1} (A simple illustrative example).  
Let 
\begin{eqnarray*}
\x & = & (x_1,\ldots,x_{n}) \in \Real^{n}, \ \Q \in \SymMat^{n} \cap \Integer^{n \times n}, \\
S & = & \left\{ \x \in \Real^{n} : \ x_1 = 1, \ x_1x_i = x_i^2 \ (i=2,\ldots,n) \right\} 
\subset \{0,1\}^{n}, \\ 
\bPhi(\x) & = & \x\x^T \in \SymMat^{n}, \ 
f(\x) \ = \ \x^T \Q \x = \inprod{\Q}{\x\x^T} = \inprod{\Q}{\bPhi(\x)}, \\ 
\widehat{S} & = & \left\{\X \in \SymMat^{n}: \X \in \SymMat^{n}_+,
 \ X_{11}=1, \   X_{1i}=X_{ii} \ (i=2,\ldots,n) \right\}.
\end{eqnarray*}
Note that $f(\x)$ with $x_1=1$ can be rewritten as
\begin{eqnarray*}
f(\x) & = & \sum_{i=2}^n\sum_{j=2}^n Q_{ij}x_ix_j + 2 \sum_{j=2}^n Q_{1j}x_j + Q_{11}. 
\end{eqnarray*} 
Thus, $\min\left\{ \x^T\Q\x : \x \in S\right\}$ corresponds to a QUBO,
and $\min\left\{ \inprod{\Q}{\X} : \X \in \widehat{S} \right\}$ to 
the standard SDP relaxation. 
Conditions (a), (b) and (e) are obviously satisfied with $\ell = n$. It is also straightforward to see 
that Conditions (d) and (f) are satisfied. 
\eexamp

\rema
Under Condition (d), problem P~\eqref{eq:NLP0} can be reformulated as 
\begin{eqnarray*}
\zeta & = & \min \left\{\inprod{\Q}{\bPhi(\x)} : \x \in S \right\}.  
\end{eqnarray*}
It is known that if $f$ is a polynomial function with integer coefficients in $\x \in \{0,1\}^n$, then we can take 
a matrix $\Q \in \SymMat^{\ell} \cap \Integer^{\ell \time \ell}$  and a mapping 
$\bPhi : \Real^n \rightarrow \SymMat^{\ell}$ 
for some $\ell$ such that Conditions (d) and (e) hold. 
Therefore, all the discussions in this section are valid for polynomial optimization problems in 
binary variables (See \cite{LASSERRE2001,LASSERRE2001b}). We mention that the authors' 
main interest is 
to develop a practical numerical method for solving large scale linearly constrained quadratic 
optimization problems in binary variables by effectively utilizing 
\em{SCIT}, 
which will be presented in Section 3.5. 
\erema

\subsection{A cutting inequality for the feasible region $S$ of  
problem P}

Let $\g \in \Integer^n $ and $\alpha \in \Integer$. 
We consider the following pair of subproblems of P:
\begin{eqnarray}
{\rm P}_{\rm a} : 
\zeta_{\rm a} & = & \min \left\{ f(\x) : \x \in S_{\rm a} \right\}, \ \mbox{where } 
S_{\rm a} = \left\{\x \in S: \g^T\x \leq \alpha \right\}, \label{eq:NLPa}  \\
{\rm P}_{\rm b} : 
\zeta_{\rm b} & = & \min \left\{ f(\x) : \x \in S_{\rm b} \right\}, \ 
\mbox{where } S_{\rm b} = \left\{\x \in S: \g^T\x \geq \alpha +1 \right\}. \label{eq:NLPb}
\end{eqnarray}
Let $\hat{\zeta}_{\rm a}$ and $\hat{\zeta}_{\rm b}$ denote the optimal values of conic 
relaxations of ${\rm P}_{\rm a}$ and ${\rm P}_{\rm b}$ such that 
$\hat{\zeta}_{\rm a} \leq \zeta_{\rm a}$ and $\hat{\zeta}_{\rm b} \leq \zeta_{\rm b}$ hold, 
respectively. More precisely, the conic relaxations of ${\rm P}_{\rm a}$ and ${\rm P}_{\rm b}$ are 
written as
\begin{eqnarray}
\widehat{{\rm P}}_{\rm a} : 
\hat{\zeta}_{\rm a} & = & 
\min \left\{ \inprod{\Q}{\X} : \X \in \widehat{S}_{\rm a} \right\}, \label{eq:ConicRelaxationa} \\
\widehat{{\rm P}}_{\rm b} : 
\hat{\zeta}_{\rm b} & = & 
\min \left\{ \inprod{\Q}{\X} : \X \in \widehat{S}_{\rm b} \right\}, \label{eq:ConicRelaxationb}
\end{eqnarray}
respectively. Here we assume that condition 
\begin{description}
\item{(${\rm f}_{\rm ab}$)} $\widehat{S}_{\rm a}$ and $\widehat{S}_{\rm b}$ are closed 
convex subsets of $\SymMat^{\ell}$ such that $\bPhi(S_{\rm a}) \subset \widehat{S}_{\rm a}$ 
and $\bPhi(S_{\rm b}) \subset \widehat{S}_{\rm b}$   
\end{description}
holds in addition to Condition (d) and (e). Hence $\hat{\zeta}_a \leq \zeta_a$ and $\hat{\zeta}_b \leq \zeta_b$. 
Since $\g^T\x^* \in \Integer$, $\x^* \in S$ lies in either $S_a$ or $S_b$. This implies 
$\zeta = \min\{\zeta_{\rm a},\zeta_{\rm b}\}$. Therefore, 
we can conclude that 
if $\eta < \hat{\zeta}_{\rm a}$, where $\eta$ denotes a known upper bound of $\zeta$  (See Condition (c)), 
then $\hat{\zeta}_{\rm b} \leq \zeta_{\rm b} = \zeta$. 
(Recall Lemma~\ref{lemma1}). In other words, problem P$_b$~\eqref{eq:NLPb} 
with the cutting inequality $\g^T\x \geq \alpha+1$ shares the same optimal value $\zeta_b = \zeta$ as the 
original problem P~\eqref{eq:NLP0}, and its conic relaxation $\widehat{\rm P}_b$~\eqref{eq:ConicRelaxationb} provides 
a lower bound $\hat{\zeta}_b$, which is at least as tight as the original lower bound $\hat{\zeta}$ for~P.

\subsection{A cutting inequality for the feasible region $\widehat{S}$ of the conic relaxation 
problem $\widehat{\rm P}$}

Let $\G \in \SymMat^{\ell} \cap \Integer^{\ell \times \ell}$ and $\alpha \in \Integer$. We consider the following pair of subproblems 
of $\widehat{\rm P}$~\eqref{eq:ConicRelaxation0}:
\begin{eqnarray}
\widetilde{\rm P}_{\rm a} : 
\tilde{\zeta}_{\rm a} & = & \min \left\{ \inprod{\Q}{\X} : \X \in \widehat{S}, \ 
\inprod{\G}{\X} \leq \alpha \right\}. \label{eq:ConicRelaxationA}\\
\widetilde{\rm P}_{\rm b} : 
\tilde{\zeta}_{\rm b} & = & \min \left\{ \inprod{\Q}{\X} : \X \in \widehat{S}, \ 
\inprod{\G}{\X} \geq \alpha + 1 \right\}. \label{eq:ConicRelaxationB}
\end{eqnarray}
Obviously, $\hat{\zeta} \leq \min \{\tilde{\zeta}_{\rm a},\tilde{\zeta}_{\rm b}\}$. 
Since $\bPhi(\x^*) \in \widehat{S} \cap \Integer^{\ell \times \ell}$ and 
$\inprod{\G}{\bPhi(\x^*)} \in \Integer$,  
we have either $\inprod{\G}{\bPhi(\x^*)} \leq \alpha$ or 
$\inprod{\G}{\bPhi(\x^*)} \geq \alpha+1$, which implies 
\begin{eqnarray*}
\tilde{\zeta}_{\rm a} \leq \inprod{\Q}{\bPhi(\x^*)} = \zeta \ \mbox{or } 
\tilde{\zeta}_{\rm b} \leq \inprod{\Q}{\bPhi(\x^*)} = \zeta. 
\end{eqnarray*}
Thus, $\min\{\tilde{\zeta}_{\rm a},\tilde{\zeta}_{\rm b}\} \leq \zeta$. Consequently,
we can conclude that 
if $\eta < \tilde{\zeta}_{\rm a}$,  then $\hat{\zeta} \leq \tilde{\zeta}_{\rm b} \leq \zeta$.  
(See Lemma~\ref{lemma2}). In other words, $\widetilde{\rm P}_b$~\eqref{eq:ConicRelaxationB} 
with the cutting inequality $\inprod{\G}{\X} \geq \alpha + 1$ provides a lower bound 
$\tilde{\zeta}_b$, 
at least as tight as the original lower bound $\hat{\zeta}$ , 
for the unknown optimal value $\zeta$ of problem P~\eqref{eq:NLP0}. 

\subsection{Successive cutting inequality technique (SCIT)}

The generation of a single cutting inequality for $S$ presented in Section~3.3 can be 
 extended in a straightforward fashion to simultaneous generation of multiple cutting inequalities.
Let $(\g_j,\alpha_j) 
\in \Integer^{n+1}$ $(j=1,\ldots,m)$. For each $j$, we consider the following subproblem of 
P~\eqref{eq:NLP0}: 
\begin{eqnarray*}
{\rm P}_{{\rm a}j} : 
\zeta_{{\rm a}j} & = & \min \left\{ f(\x) : \x \in S_{{\rm a}j}\right\}, \ 
\mbox{where } S_{{\rm a}j} = \left\{ \x \in S: 
\g_j^T\x \leq \alpha_j \right\},  
\end{eqnarray*}
and a conic relaxation of ${\rm P}_{{\rm a}j}$
\begin{eqnarray}
\widehat{{\rm P}}_{{\rm a}j} : 
\widehat{\zeta}_{{\rm a}j} & = & \min \left\{ \inprod{\Q}{\X} : \X \in \widehat{S}_{{\rm a}j}
\right\}.  \label{eq:ConicRelaxationaj}
\end{eqnarray}
Let $J = \{ j : \eta < \hat{\zeta}_{{\rm a}j} \}$. Then, 
by adding the cutting inequalities $\g_j^T\x \geq \alpha_j+1$ $(j \in J)$ to problem P, 
we obtain 
\begin{eqnarray}
\mbox{P$^1$: } \zeta & = & \min \left\{f(\x) : \x \in S^1 \right\},  \label{eq:NLP1} 
\end{eqnarray}
where 
$ 
S^1 = \left\{ \x \in S : \g_j^T\x \geq \alpha_j+1 \ (j \in J) \right\}, 
$ 
and its conic relaxation:
\begin{eqnarray}
\widehat{\rm P}^1:  \hat{\zeta}^1 & = & \min \left\{\inprod{\Q}{\X} : \X \in \widehat{S}^1 \right\}.   \label{eq:ConicRelaxationOne} 
\end{eqnarray}
Problem ${\rm P}^1$~\eqref{eq:NLP1} shares the same optimal value $\zeta$ 
and optimal solution $\x^*$ with 
problem ${\rm P}$, and $\hat{\zeta} \leq \hat{\zeta}^1 \leq \zeta$.  
Note that problems $\widehat{\rm P}_{{\rm a}j}$ $(j=1,\ldots,m)$ 
can be solved independently in parallel.

The simultaneous generation of  multiple cutting inequalities above   
can be applied  now
to the feasible region $S^1$ of 
problem P$^1$, and a new problem can be constructed as 
\begin{eqnarray*}
{\rm P}^2 : \zeta & = & \left\{ f(\x) : \x \in S^2 \right\},   
\end{eqnarray*}
which is equivalent to P and P$^1$, and its conic relaxation problem 
\begin{eqnarray*}
\widehat{\rm P}^2 : \hat{\zeta}^2 & = & \left\{ \inprod{\Q}{\X} : \X \in \widehat{S}^2 \right\}.  
\end{eqnarray*}
We continue this process successively to generate 
a sequence of conic relaxation problems 
$\{\widehat{\rm P}^k : k=1,2,\ldots\}$ (of $\{{\rm P}^k : k=1,2,\ldots\}$) and a sequence 
$\{\hat{\zeta}^k \ (k=1,\ldots,\ ) \}$ of their optimal values, which serve as 
lower bounds for $\zeta$ such that 
$\hat{\zeta}\leq \hat{\zeta}^k \leq \hat{\zeta}^{k+1} \leq \zeta$ $(k=0,1,\ldots )$. 
This entire process constitutes   {\em SCIT}.

For a similar extension of generating a single cutting inequality for $\widehat{S}$  presented in Section~3.4 
to SCIT, let $\G_j \in \SymMat^{\ell} \cap \Integer^{\ell \times \ell}$ and $\alpha_j \in \Integer$ 
$(j=1,\ldots,m)$. We consider 
\begin{eqnarray*}
\widetilde{\rm P}_{{\rm a}j} : \tilde{\zeta}_{{\rm a}j} &=& 
\min\left\{\inprod{\Q}{\X}: \X \in \widehat{S}, \ 
\inprod{\G_j}{\X} \leq \alpha_j\right\}
\end{eqnarray*}
for each $j=1,\ldots,m$. Let $J =\{j : \eta < \tilde{\zeta}_{{\rm a}j} \}$. Then, 
the cutting inequalities $\inprod{\G_j}{\X} \geq \alpha_j + 1$ 
$(j \in J)$ can be applied to problem $\widehat{\rm P}$~\eqref{eq:NLP0}: 
\begin{eqnarray*}
\widetilde{\rm P}^1 : \tilde{\zeta}^1 & = & \min\left\{ \inprod{\Q}{\X}: \X \in \widetilde{S}^1 \right\},
\end{eqnarray*}
where $\widetilde{S}^1 = \left\{\X \in \widehat{S} : \inprod{\G_j}{\X} \geq \alpha_j + 1 \ 
(j \in J) \right\}$. 
As a result, we obtain that $\hat{\zeta} \leq \tilde{\zeta}^1 \leq \zeta$. We note that conic relaxation 
problems $\widetilde{\rm P}_{{\rm a}j}$ $(j=1,\ldots,m)$ can be solved independently in parallel. 

Now, the discussion above is applied to  the feasible region $\widetilde{S}^1$ of 
problem $\widetilde{\rm P}^1$, and a new conic relaxation problem 
$\widetilde{\rm P}^2$ of problem P  is constructed as:
\begin{eqnarray*}
\widetilde{\rm P}^2 : \tilde{\zeta}^2 & = & \min\left\{ \inprod{\Q}{\X}: \X \in \widetilde{S}^2 \right\}. 
\end{eqnarray*}
Continuing this process, a sequence of conic relaxation problems 
$\left\{\widetilde{\rm P}^1,\widetilde{\rm P}^2,\cdots\right\}$ and a sequence 
$\{\widetilde{\zeta}^1,\widetilde{\zeta}^2,\ldots \}$ 
of their optimal values are generated such that 
$\hat{\zeta} \leq \widetilde{\zeta}^k \leq \widetilde{\zeta}^{k+1} \leq \zeta$ $(k=1,2,\ldots)$.  
Thus generating a single cutting inequality for $\widehat{S}$ has been extended to 
SCIT for $\widehat{S}$.

\section{An application to quadratic unconstrained binary optimization problem (QUBO)}

We demonstrate in this section that 
SCIT 
presented in Section 3.5 has  promising prospects 
for solving 
QUBOs. More precisely, we show that the exact optimal values of some  QUBO instances can be obtained
by simply applying SCIT, without a well-designed numerical method for solving QUBOs. 

It should be mentioned that SCIT 
needs to be eventually combined with other practical numerical methods, such as the branch-and-bound 
method \cite{GLOVER1998,GUIMARAES2020,KRISLOCK2014,KRISLOCK2017},  heuristic methods including the tabu search \cite{GLOVER1998} and 
the genetic algorithm \cite{MITCHELL1998} to
solve QUBOs and other 
binary QOPs.
Before designing a specific numerical method using SCIT
and conducting extensive numerical experiments on a parallel machine with a large number of cores, 
we investigate the performance of SCIT 
on $60$ QUBO instances with dimensions 
100 - 250 from BIQMAC \cite{BIQMAC}. Obviously,  a lot of flexibility exists in implementing SCIT  
and many details 
should be determined. In the subsequent discussion, we choose some specific values 
for SCIT 
to just carry out numerical experiments. Those settings are not to propose a numerical method for solving QUBOs. 
Nevertheless, 
 the exact  optimal values could be attained for 70\% of the 
$60$ QUBO instances in $10$ iterations (see Section 4.4). 

\subsection{A QUBO}

We consider a 
QUBO:
\begin{eqnarray}
\zeta &=& \min \left\{\u^T\R\u : \u \in \{0,1\}^m \right\},   \label{eq:QUBO00}
\end{eqnarray}
where $\R \in \SymMat^m$. Introducing a slack variable 
vector $\v \in \{0,1\}^m$, we transform the QUBO to 
\begin{eqnarray}
\mbox{P: } \zeta & = & \min \left\{\u^T\R\u : \x = (\u,\v) \in S \right\}, 
\label{eq:QUBO0} 
\end{eqnarray}
where $S = \left\{ \x = (\u,\v) \in \{0,1\}^{2m} : \u + \v = \e \right\}$ and 
$\e$ denotes the $m$-dimensional column vector of $1$'s. 
It is known that 
introducing the slack variabl vector $\v \in \{0,1\}^m$ is crucial to strengthen the 
conic relaxation of QUBO~\eqref{eq:QUBO0} (see, for example, \cite[Section 6.1]{ITO2017}). 
We note that $(x_1,\ldots,x_m)$ corresponds to 
$\u \in \Real^m$ and $(x_{m+1},\ldots,x_{2m})$ to $\v \in \Real^m$, and that 
$\x = (\u,\v) \in S$ implies the complementarity 
$x_i x_{m+i} = u_iv_i = 0$ between $x_i = u_i$ and $x_{m+i} = v_i$ $(i=1,\ldots,m)$; hence
 $\sum_{i=1}^{2m} x_i = m$ holds for every $\x = (\u,\v) \in S$. 
 Let $\x^* = (\u^*,\v^*)$ be an unknown optimal 
solution of QUBO~\eqref{eq:QUBO0}, and $\eta$ a known upper bound of the optimal 
value $\zeta$. 

\subsection{Cutting inequalities in the number of $1$'s in $(\u,\v) \in \{0,1\}^{2m}$}

For each $I \subset \{1,\ldots,2m\}$ and $\alpha \in \{0,1,\ldots,|I|-1\}$ where $|I|$ denotes 
the number of elements of $I$, we consider the following type of cutting inequality 
for the feasible region $S$ of QUBO~\eqref{eq:QUBO0}:
\begin{eqnarray*}
\sum_{i \in I} x_i \geq \alpha+1, 
\end{eqnarray*}
which together with $\x = (\u,\v) \in S$ 
requires that the number of $1$'s among $x_i$ $(i \in I)$ 
is at least $\alpha+1$. 
In particular, if we take $I = \{i\}$ with $i \in \{1,\ldots,m\}$ (or $i \in \{m+1,\ldots,2m\}$) 
and $\alpha + 1 = 1$, the cutting inequality $\sum_{i \in I} x_i \geq \alpha + 1$ requires $u_i=1$ 
and $v_i=0$ (or $u_i=0$ and $v_i=1$). Thus, it is possible to fix $u_i$ to $1$ (or $0$)  and reduce 
the size of QUBO if the inequality is shown to be valid for $\x = \x^*$. 

\subsection{An experimental method using SCIT}

To initialize the sequence 
\begin{eqnarray}
\left\{ \{(I,\alpha^k_I,\beta^k_I) : I\in\IC^k\} : k=0,1,\ldots, \right\}, \label{eq:sequence0}
\end{eqnarray}
which is to be generated,  set 
\begin{eqnarray*}
\IC^0 & = & \mbox{a family of nonempty subsets of $\{1,\ldots,2m\}$}, \\  
\alpha^0_I & = & 0 \ \mbox{ and } \beta^0_I \ = \ \lfloor \gamma |I| \rfloor \ 
\mbox{for every $I \in \IC^0$}, 
\end{eqnarray*}
where $\gamma = 1/3$ is used in the preliminary numerical experiment reported  in Section~4.4.
For each $k=0,1,\ldots$, we consider 
\begin{eqnarray}
{\rm P}^k : \zeta^k & = & \min \left\{\u^T\R\u : 
\x = (\u,\v) \in S, \ 
\displaystyle \sum_{i \in I} x_i \geq \alpha^k_I \ (I \in \IC^k) 
\right\}, \label{eq:QUBOk}
\end{eqnarray}
and its conic relaxation problem $\widehat{\rm P}^k$ with the optimal value $\hat{\zeta}^k$. 

Let $k=0$. Since $\alpha^k_I = 0$ for every $I \in \IC^k$,  the inequalities  
$\sum_{i \in I} x_i \geq \alpha^k_I \ (I \in \IC^k)$ obviously hold for every $\x = (\u,\v) \in S$ 
and  problem P$^k$~\eqref{eq:QUBOk} is equivalent to problem P~\eqref{eq:QUBO0}. Hence, 
\begin{eqnarray}
\mbox{$\x^*$ remains an optimal solution of ${\rm P}^k$} \ \mbox{ and } 
\hat{\zeta}^{k} \leq \zeta^k = \zeta. 
\label{eq:iterationk1}
\end{eqnarray}
Assuming that~\eqref{eq:iterationk1} holds for some iteration $k \in \{0,1,\ldots\}$, 
we show how to update 
$ \{(I,\alpha^k_I,\beta^k_I) : I\in\IC^k\} $ 
 to 
$ \{(I,\alpha^{k+1}_I,\beta^{k+1}_I) : I\in\IC^{k+1}\}$  
so that 
\begin{eqnarray}
\mbox{$\x^*$ remains an optimal solution of ${\rm P}^{k+1}$} \ \mbox{ and } 
\hat{\zeta}^{k} \leq \hat{\zeta}^{k+1} \leq \zeta^{k+1} = \zeta. 
\label{eq:iterationkplus1}
\end{eqnarray}
In the numerical experiment whose results are reported in Section 4.4, 
the 
Lagrangian-DNN relaxation \cite{KIM2013} (see also \cite{ARIMA2017,ARIMA2018}) 
for $\widehat{\rm P}^k$ and $\widehat{\rm P}_a^k(I')$ 
described below was employed, and  NewtBracket \cite{KIM2020} (the Newton-bracketing method \cite{KIM2019a}) 
was applied to them for their optimal values $\hat{\zeta}^k$ and $\hat{\zeta}^k(I')$, respectively. 

For simplicity of discussion, we first deal with the case where $\IC^{k+1} = \IC^{k}$. 
For each $I' \in \IC^{k}$, we consider the following problem: 
\begin{eqnarray*}
{\rm P}^k_a(I') : \zeta^k(I') & = & \min \left\{\u^T\R\u : 
\begin{array}{l}
\x = (\u,\v) \in S, \ 
\displaystyle \sum_{i \in I} x_i \geq \alpha^k_I \ (I \in \IC^{k} \backslash I'),\\ [3pt]
\displaystyle \sum_{i \in I'} x_i \leq \alpha^k_{I'} + \beta^k_{I'} 
\end{array} 
\right\}, 
\end{eqnarray*}
and solve its conic relaxation $\widehat{\rm P}^k_a(I')$ to compute its optimal value 
$\hat{\zeta}^k_a(I')$. If $\eta < \hat{\zeta}^k_a(I')$, let  
\begin{eqnarray*}
\alpha^{k+1}_{I'} =  \alpha^{k}_{I'} + \beta^k_{I'} + 1 \ \mbox{ and } 
\beta^{k+1}_{I'}   =  \min \left\{ \beta^{k}_{I'}, |I'| - \alpha^{k+1}_{I'} - 1 \right\}.  
\end{eqnarray*}
Otherwise, let 
\begin{eqnarray*}
\alpha^{k+1}_{I'} =  \alpha^{k}_{I'} \ \mbox{and } 
\beta^{k+1}_{I'}  =  \lfloor \beta^{k}_{I'}/2 \rfloor.  
\end{eqnarray*}
Thus $ \{(I,\alpha^k_I,\beta^k_I) : I\in\IC^k\} $ has been updated to 
$ \{(I,\alpha^{k+1}_I,\beta^{k+1}_I) : I\in\IC^{k+1}\}$.  
From the discussion in Sections~3.3 and~3.5, we  see that~\eqref{eq:iterationkplus1} holds. 

\subsubsection{An algorithm for generating $\IC^0$} 

If $\x^*$ were known in advance, it would be easy to construct an ideal cutting inequality of the 
form $\sum_{i \in I}x_i \geq \alpha$ such that 
$
\left\{ \x =(\u,\v) \in S :  \sum_{i \in I}x_i \geq \alpha \right\} = \{\x^*\}. 
$
In fact, we could take $I = \{ i : x^*_i = 1\}$ and $\alpha = m$, which  is 
impossible. If a branch-and-bound method, for instance, is applied to solve QUBO~\eqref{eq:QUBO0}, 
then more accurate information on the location of optimal solutions becomes available as it proceeds.
In such a case, it is reasonable to incorporate such information into 
$\{ (I, \alpha^k_I,\beta^k_i) :  I \in \IC^k \}$. 
This will be discussed in Section~4.3.2. 

For the case where no information on 
the location of the optimal solutions of QUBO~\eqref{eq:QUBO0} is available, we propose `to 
distribute the cutting inequalities uniformly'.   
There  still remains a great deal of flexibility in  choosing $\IC^0$ to initialize the sequence~\eqref{eq:sequence0}.
In general, as the members of $\IC^0$ increase,  a tighter lower bound $\hat{\zeta}^k$ 
for the optimal value $\zeta$ of P~\eqref{eq:QUBO0} at each $k$th iteration can be expected. 

Let us show a simple example of $\IC^0$ below, which may provide an idea for a  general 
choice of~ $\IC^0$. 
\begin{description}
\item{Step 0: } Let $r=0$. $\JC^0 = \left\{\{1,\ldots,m\}\right\}$. 
\item{Step 1: } If $|J| = 1$ for all $J \in \JC^r$ then let
\begin{eqnarray*}
\JC & = & \bigcup_{p=0}^r \JC^p, \ 
\IC^0 \ = \ \JC \bigcup \left\{\{m+j: j \in J\} : J \in \JC \right\},  
\end{eqnarray*}
and stop. 
\item{Step 2: } Let $\JC^{r+1} = \emptyset$. For every $J \in \JC^r$ with $|J| \geq 2$, 
choose two subsets $J_1$ and $J_2$ of $J$ (randomly) such that $J_1 \cup J_2 = J$ and 
$|J_1| = |J_2| = \lceil |J|/2 \rceil$ (then $|J_1 \cap J_2| \leq 1$),  and add them 
to $\JC^{r+1}$. 
\item{Step 3: } Let $r = r+1$ and go to Step 1.
\end{description}

If $n=4$, the above algorithm generates 
\begin{eqnarray*}
& & \JC^0 = \left\{\{1,2,3,4\}\right\}, \ \JC^1 = \left\{\{1,3\},\{2,4\}\right\},\ 
\JC^2 = \left\{\{1\},\{2\},\{3\},\{4\}\} \right\}, \\ 
& & \JC = \bigcup_{p=1}^2 \JC^p = \left\{\{1,2,3,4\},\{1,3\},\{2,4\},\{1\},\{2\},\{3\},\{4\} \right\}, \\
& & \IC^0 = \JC \bigcup \left\{\{5,6,7,8\},\{5,7\},\{6,8\},\{5\},\{6\},\{7\},\{8\}\right\}.   
\end{eqnarray*}

\subsubsection{Adding more cutting  inequalities at each iteration}

We now consider the case where some information of the location of the optimal solutions of 
P~\eqref{eq:QUBO0} is available, and discuss  
how it can be used in the construction of $\IC^{k+1}$.  
Suppose that we have generated $\{(I,\alpha^k_I,\beta^k_I): I \in \IC^k\}$ 
at which~\eqref{eq:iterationk1} holds. Assume that the information is given as 
an $\bar{\x} = (\bar{\u},\bar{\v}) \in [0,1]^{2m}$ but not necessary $\bar{\x} \in \{0,1\}^{2m}$, which 
is obtained from an optimal solution of a conic (SDP and DNN) relaxation of 
P$^k$~\eqref{eq:QUBOk}. We note that $\bar{\x} = (\bar{\u},\bar{\v})$ 
satisfies $\bar{\u}+\bar{\v}=\e$ approximately,  but  may not satisfy the complementarity 
$\bar{u}_i\bar{v}_i = 0$ $(i=1,\ldots,m)$. 
In this case, for the computation of an approximate solution $\hat{\u} \in \{0,1\}^{m}$ 
of QUBO~\eqref{eq:QUBO00}, rounding is frequently applied  to $\bar{\u} \in [0,1]^{m}$ and/or 
a heuristic method such as the tabu search \cite{GLOVER1998} and 
the genetic algorithm \cite{MITCHELL1998} to 
QUBO~\eqref{eq:QUBO00} with the initial solution $\bar{\u} \in [0,1]^{m}$.

For the construction of a family $\IC^+$ of subsets of $\{1,\ldots,2m\}$ 
to be added to $\IC^{k}$, 
each $\bar{u}_i \in [0,1]$ and $\bar{v}_i \in [0,1]$ are regarded to represent the probability 
Pr$\{ u^*_i = 1 \}$  and Pr$\{ v^*_i = 1 \}$), respectively,  for the unknown 
optimal solution $\x^* = (\u^*,\v^*)$ of P~\eqref{eq:QUBO0}, and 
$2q$ points $\u^p \in \{0,1\}^s$  $(p=1,\ldots,q))$ and  $\v^p \in \{0,1\}^s$ $(p=1,\ldots,q))$ 
are generated randomly  using the probability. 
Then, let 
\begin{eqnarray*}
\IC^+ & = & \left\{ \{ i : u^p_i = 1 \} : p=1,\ldots,q \right\} \bigcup 
 \left\{ \{ i +m : v^p_i = 1 \} : p=1,\ldots,q \right\}, \\ 
\IC^{k+1} & = & \IC^k \bigcup \IC^+.  
\end{eqnarray*}
We took $q = 10$ in the numerical experiment presented in Section 4.4. 

Now, we consider the case where an approximate optimal solution 
$\hat{\x} = (\hat{\u},\hat{\v}) \in S$ 
of P~\eqref{eq:QUBO0}, 
 which is likely to be optimal but has not been proved to be optimal, is known   
 with the objective value $\eta = \bar{\u}^T\R\bar{\u}$.
 Note that $\zeta \leq \eta$ is guaranteed. 
Such a case frequently occurs when we try to solve P~\eqref{eq:QUBO0} by a high performance 
heuristic method. Let $\widehat{I} = \{i : \hat{x}_i = 1\}$. Then $|\widehat{I}| = m$.
 For every nonempty subset $J$ of 
$\widehat{I}$ and $\alpha_{J} \in \{0,\ldots,|J|-1\}$, consider the following problem: 
\begin{eqnarray*}
{\rm P}^{k}_a(J) : \zeta^{k}_{\rm a}(J) & = & \min \left\{\u^T\R\u : 
\begin{array}{l}
\x = (\u,\v) \in S, \ 
\displaystyle \sum_{i \in I} x_i \geq \alpha^k_I \ (I \in \IC^{k}),\\ [3pt]
\displaystyle \sum_{i \in J} x_i \leq \alpha_{J} 
\end{array} 
\right\}, 
\end{eqnarray*}
and solve its conic relaxation $\widehat{\rm P}^k_a(J)$ to compute its optimal value 
$\hat{\zeta}^k_a(J)$. If $\eta < \hat{\zeta}^k_a(J)$ holds, then we know that 
$\sum_{i \in J} x_i \geq \alpha_{J}+1$ is a cutting inequality for the set of optimal solutions 
of P~\eqref{eq:QUBO0}. Moreover, if $\alpha_{J}+1 = |J|$, 
 then  $x_i$ can be fixed to $x_i = 1$
for all $i \in J$. 
If we take $J = \widehat{I}$ and $\alpha_{J}+1 =|\widehat{I}| =  m$, then 
$\eta < \hat{\zeta}^k_a(\widehat{I})$ 
provides a certificate for $\hat{\x} = (\hat{\u},\hat{\v}) \in S$ to be the unique optimal solution 
of P~\eqref{eq:QUBO0}. Therefore, it is reasonable to include  $\widehat{I}$  and/or 
some of its subsets $J$ in $\IC^+$.

\rema \label{remark:branching} 
For the case above,  
a branching  can be used instead of cutting inequalities to 
efficiently solve P~\eqref{eq:QUBO0} to optimality. 
More precisely,
for each 
$\alpha = 0,1,\ldots,m$, let  
\begin{eqnarray*}
S(\alpha) = \left\{ \x =(\u,\v) \in S : \sum_{i \in \widehat{I}} x_i = \alpha \right\}.
\end{eqnarray*}
Then,  P~\eqref{eq:QUBO0} is branched into $1+m$ subproblems 
\begin{eqnarray*}
       {\rm P}(\alpha) : \zeta(\alpha) = \{\u^T\R\u : \x = (\u,\v) \in S(\alpha) \} \ 
       (\alpha = 0,1,\ldots,m). 
\end{eqnarray*}
The important  features of this branching are:\vspace{-2mm} 
\begin{itemize}
\item As $\alpha$ increases from $0$ to $m$, $|S(\alpha)|$ deceases to $|S(m)| = 1$, 
so that subproblem ${\rm P}(\alpha)$ with a larger $\alpha$ is easier to solve. \vspace{-2mm} 
\item As $\alpha$ decreases from $m$ to $0$, the optimal value $\zeta(\alpha)$ 
of subproblem ${\rm P}(\alpha)$ is expected to increase, larger than $\eta$. 
As a result, the possibility that ${\rm P}(\alpha)$ is pruned by the lower bounding procedure 
 using its conic relaxation is increased. \vspace{-2mm} 
\end{itemize}
Further investigation of this branching is beyond the scope of the paper. 
It will be investigated in our future work. 
\erema

\subsection{Preliminary numerical results}

An experimental method on QUBOs has been described in Section 4.3 
for evaluating the performance of 
SCIT 
presented in Section~3.5.
We applied 
the method to $60$ QUBO instances from BIQMAC \cite{BIQMAC}:
\begin{eqnarray*}
& & \mbox{bqp100-1,$\ldots$,bqp100-10,be120.3.1,$\ldots,$be120.3.10,be120.8.1,$\ldots,$be120.8.10},\\
& & \mbox{be150.3.1,$\ldots,$be150.3.10,be150.8.1,$\ldots,$be150.8.10,bqp250-1,$\ldots$,bqp250-10.}
\end{eqnarray*}
The experiments were performed on iMac Pro with Intel Xeon W CPU (3.2 GHZ), 8 cores and 128 GB memory for the instances with dimensions 100, 120 and 150, 
and Intel Xeon 4216 2 CPUs with 32 cores and 128 GB memory for the instances 
with dimension 250. 

As the optimal value $\zeta$ of each instance above is known, 
 its upper bound $\eta$  was set to $\zeta$ 
to ensure the best performance of SCIT. 
Recall that $\eta$ is used in the certificate 
$\eta < \hat{\zeta}_{{\rm a}j}$  for $\g_j^T\x \geq \alpha_j+1$ to be a cutting plane, where 
$\hat{\zeta}_{{\rm a}j}$ 
denotes the optimal value of $\widehat{\rm P}_{{\rm a}j}$~\eqref{eq:ConicRelaxationaj}. 
As a smaller $\eta \geq \zeta$ is chosen
(or  $\eta$ closer to $\zeta$), more inequalities can become valid cutting inequalities. 
Thus, $\zeta$
is the best choice in our experiment.
In the case where $\zeta$ is not known, $\eta$ is usually obtained by a heuristic method
as it needs to be the best known upper bound of the optimal value $\zeta.$

For the conic relaxation 
$\widehat{\rm P}^k$ of P$^k$ and $\widehat{\rm P}^k_{\rm a}(I')$ of 
${\rm P}^k_{\rm a}(I')$ $(I' \in \IC^k)$, 
the Lagrangian-DNN relaxation \cite{KIM2013,ARIMA2017,ARIMA2018} of P$^k$ and
P$^k_{\rm a}(I')$ $(I' \in \IC^k)$ was employed, respectively, and NewtBracket \cite{KIM2020} (the Newton-bracketing method \cite{KIM2019a}) as a numerical method to compute their optimal values 
$\hat{\zeta}^k$ 
and $\hat{\zeta}^k_{\rm a}(I')$ $(I' \in \IC^k)$. $\IC^0$ was constructed as described in 
Section 4.3.1. 
An approximate optimal solution 
$\overline{\X}$ of $\widehat{\rm P}^k$ was also computed, and an $\bar{\x} \in [0,1]^{2m}$ 
described in Section 4.3.2 from $\overline{\X}$ was obtained for the information on the location 
of the optimal solutions P$^k$. 
We used $q = 10$ for the additional cut inequalities associated with $\IC^+$.
$\IC^0$ contains about $4m$ subsets of $\{1,\ldots,2m\}$, so that approximately $4m$ valid 
inequalities of the form $\sum_{i \in I} x_i \geq \alpha$ were prepared prior to the $0$th 
iteration, and $|I^+| = 10$ inequalities were added prior to  the $k$th iteration $(k=1,2,\ldots)$.  

Each $k$th iteration consists of two phases: the first one for solving $\widehat{\rm P}^k$ and 
the second one for solving $\widehat{\rm P}^k_a(I')$ $(I' \in \IC^k)$. After solving the second
problems, 
$ \{(I,\alpha^{k}_I,\beta^{k}_I) : I\in\IC^{k}\}$ was updated to 
$ \{(I,\alpha^{k+1}_I,\beta^{k+1}_I) : I\in\IC^{k+1}\}$ as described in Sections 4.3. 
Since $\IC^0$ contains a singleton $I = \{i\}$ $(i=1,\ldots,2m)$, 
the variable $x_i$ was fixed to $1$ when the inequality $x_i \geq 1$ became 
valid cutting inequality, and the number of free variables was reduced 
among $u_1=x_1,\ldots,u_m=x_m$ (= the 
dimension of subQUBO denoted by $d^k$ in Table 1) as well as the redundant 
cutting inequalities were removed. 
 
The iteration terminated when $\hat{\zeta}^k$ attained the optimal value $\zeta$ in $10$ iterations 
$k=0,1,\ldots,9$ or $k$ reached $9$. 
Among $60$ QUBO instances, $42$ cases attained the optimal value within $10$ iterations, 
{\it i.e.}, $\hat{\zeta}^k = \zeta$ for some $k <= 9$. 
The  other $18$ cases failed to obtain the exact optimal value.
Table 1 shows the numerical results on $42$ successful instances.

The following two aspects are crucial for evaluating the performance 
of SCIT: 
\vspace{-2mm}
\begin{description}
\item{(i) } 
How many 
variables 
$u_1=x_1,\ldots,u_m=x_m$ are fixed to either $0$ or $1$, which can be measured by 
the decrease of $d_k$ $(k=0,1,\ldots)$.\vspace{-2mm} 
\item{(ii) } Improvement in the lower bound of optimal value $\zeta$, which can be observed 
by the increase of 
$\hat{\zeta}^k$  $(k=0,1,\ldots)$.\vspace{-2mm} 
\end{description}
Overall, the method worked effectively in terms of the two aspects, but less effectively for larger 
dimensional cases; it took more iterations to attain a smaller $d^k$ and a tighter 
$\hat{\zeta}^k$ to $\zeta$. In practice, (i) is an important aspect of SCIT
when it is combined with a 
numerical method for solving QUBOs. For the QUBO instance bqp250-2 in our numerical experiment,  
  $237 (=250-13)$ variables among $u_1,\ldots,u_{250}$ 
of QUBO~\eqref{eq:QUBO0} were fixed to $0$ or $1$ after $10$ iterations, 
so the resulting subQUBO 
with $13$ variables was
easy to solve.

The $18$ instances where  $\hat{\zeta}^k$ could not attain 
$\zeta$ in $10$ iterations (or $k \leq 9$) 
are not included in Table 2. 
As mentioned before, SCIT
alone cannot be 
a numerical method for solving QUBOs. To successfully solve the 18 instances, a numerical method combining
 SCIT 
 with other methods should  be implemented.

 We want to highlight that the exact optimal values of the 42 QUBO instances, 
 $70$ \% of the $60$ instances to which the method was applied,
 could be obtained in the numerical experiments.
 The numerical results reported here,  though  limited, 
present the promising potential of SCIT, especially when it is combined with 
the branch-and-bound method \cite{GLOVER1998,GUIMARAES2020,KRISLOCK2014,KRISLOCK2017} and 
heuristic methods \cite{GLOVER1998,MITCHELL1998} for solving QUBOs.

\begin{table}[htp]
\begin{center} 
\caption{Numerical results on SCIT applied to 
42 QUBO instances from BIQMAC \cite{BIQMAC}. 
}
\label{table:numRes1}
\vspace{2mm} \hspace*{-0.5cm}
\scalebox{0.60}{
\begin{tabular}{|l|r|r|r|r|r|r|r|r|r|r|r|r|r|}
\hline 
                  &                  & 
\multicolumn{10}{|c|}{$d^k$  (Dim. of subQUBO ${\rm P}^k$ whose L-DNN relaxation 
$\widehat{\rm P}^k$ to be solved), \ 
$\hat{\zeta}^k$ (the optimal value of $\widehat{\rm P}^k$)}  \\
QUBO   &  Opt.Val   & $k=0$                    & $k=1$               & $k=2$  & $k=3$  & $k=4$  & $k=5$ & $k=6$ & 
$k=7$  & $k=8$ & $k=9$ \\
\hline
bqp100-1  & -7970 & 100, -8036 & 42, -7970 &   &   &   &   &   &   &   &   \\
bqp100-2 & -11036 & 100, -11036 &                   &    &        &       &       &       &     &       &   \\ 
bqp100-3 & -12723 & 100, -12723 & & & & & & & & &   \\
bqp100-4 & -10368 & 100, -10368 & & & & & & & & &    \\
bqp100-5 & -9083 & 100, -9083 & & & & & & & & &    \\
bqp100-6  & -10210 & 100, -10341 & 56, -10291 & 38, -10270 & 27, -10248 & 17, -10220 & 12, -10210 &   &   &   &   \\
bqp100-7  & -10125 & 100, -10159 & 36, -10125 &   &   &   &   &   &   &   &   \\
bqp100-8 & -11435 & 100, -11435 & & & & & & & & &    \\
bqp100-9 & -11455 & 100, -11455 & & & & & & & & &    \\
bqp100-10 & -12565 & 100, -12565 & & & & & & & & &    \\
\hline 
be120.3.1  & -13067 & 120, -13343 & 93, -13268 & 78, -13204 & 58, -13135 & 35, -13075 & 9, -13067 &   &   &   &   \\
be120.3.2  & -13046 & 120, -13163 & 44, -13046 &   &   &   &   &   &   &   &   \\
be120.3.3  & -12418 & 120, -12609 & 75, -12477 & 30, -12418 &   &   &   &   &   &   &   \\
be120.3.4  & -13867 & 120, -14039 & 71, -13939 & 31, -13868 & 2, -13867 &   &   &   &   &   &   \\
be120.3.5  & -11403 & 120, -11558 & 59, -11407 & 12, -11403 &   &   &   &   &   &   &   \\
be120.3.6  & -12915 & 120, -13022 & 46, -12915 &   &   &   &   &   &   &   &   \\
be120.3.7  & -14068 & 120, -14128 & 27, -14068 &   &   &   &   &   &   &   &   \\ 
be120.3.8  & -14701 & 120, -14812 & 40, -14701 &   &   &   &   &   &   &   &   \\
be120.3.10  & -12201 & 120, -12413 & 83, -12298 & 46, -12202 & 4, -12201 &   &   &   &   &   &   \\
\hline
be120.8.2  & -18827 & 120, -19351 & 102, -19271 & 96, -19167 & 83, -19065 & 61, -18903 & 29, -18827 &   &   &   &    \\
be120.8.3  & -19302 & 120, -19791 & 102, -19653 & 80, -19509 & 50, -19396 & 19, -19302 &   &   &   &   &   \\
be120.8.4  & -20765 & 120, -21063 & 65, -20824 & 12, -20765 &   &   &   &   &   &   &   \\
be120.8.5  & -20417 & 120, -20677 & 46, -20457 & 21, -20417 &   &   &   &   &   &   &   \\
be120.8.6  & -18482 & 120, -18954 & 98, -18804 & 74, -18615 & 35, -18482 &   &   &   &   &   &   \\
be120.8.9  & -18195 & 120, -18685 & 101, -18539 & 80, -18384 & 49, -18231 & 23, -18195 &   &   &   &   &   \\
be120.8.10  & -19049 & 120, -19380 & 66, -19157 & 28, -19055 & 6, -19049 &   &   &   &   &   &   \\
\hline
be150.3.1  & -18889 & 150, -19202 & 117, -19098 & 89, -18978 & 50, -18889 &   &   &   &   &   &    \\ 
be150.3.2  & -17816 & 150, -18200 & 129, -18123 & 110, -18069 & 102, -17982 & 76, -17861 & 36, -17816 &   &   &   &     \\
be150.3.3  & -17314 & 150, -17510 & 79, -17315 & 13, -17314 &   &   &   &   &   &   &   \\ %
be150.3.4  & -19884 & 150, -20080 & 82, -19917 & 27, -19884 &   &   &   &   &   &   &   \\ 
be150.3.5  & -16817 & 150, -17216 & -139, 17159 & 123, -17092 & 104, -16998 & 85, -16930 & 63, -16846 & 46, -16817 &   &   &   \\
be150.3.7  & -18001 & 150, -18385 & 132, -18331 & 117, -18248 & 93, -18141 & 62, -18052 & 27, -18001 &   &   &   &   \\
\hline
be150.8.4  & -26911 & 150, -27685 & 142, -27611 & 136, -27561 & 133, -27476 & 130, -27455 & 128, -27423 & 122, -27371 & 113, -27257 & 92, -27108 & 58, -26911 \\
be150.8.5  & -28017 & 150, -28634 & 116, -28470 & 96, -28343 & 68, -28198 & 44, -28076 & 26, -28022 & 6, -28017 &   &   &   \\
be150.8.10  & -28374 & 150, -29125 & 139, -29071 & 132, -28950 & 118, -28894 & 112, -28810 & 102, -28691 & 84, -28607 
& 69, -28492 & 43, -28374 &   \\ 
\hline
bqp250-1  & -45607 & 250, -46244 & 214, -46102 & 181, -45927 & 134, -45719 & 52, -45607 &   &   &   &   &   \\
bqp250-2  & -44810 & 250, -45585 & 241, -45551 & 230, -45469 & 219, -45413 & 211, -45346 & 199, -45281 & 178, -45186 & 156, -45085 & 109, -44841 & 13, -44810 \\
bqp250-3  & -49037 & 250, -49457 & 163, -49144 & 56, -49037 &   &   &   &   &   &   &  \\
bqp250-4  & -41274 & 250, -42009 & 237, -41909 & 223, -41805 & 198, -41668 & 166, -41478 & 92, -41274 &   &   &   &  \\
bqp250-5  & -47961 & 250, -48431 & 164, -48153 & 82, -48040 & 45, -47961 &   &   &   &   &   &   \\
bqp250-7  & -46757 & 250, -47378 & 215, -47236 & 182, -47086 & 130, -46904 & 59, -46757 &   &   &   &   &   \\
\hline
\end{tabular}
}
\end{center}
\end{table}

\section{Concluding Remarks}


We have presented SCIT, a very flexible framework, to generate 
effective cutting inequalities 
for strengthening conic relaxations   
for computing lower bounds of the optimal value of a binary QOP. 
%
%
%
%
%
To be able to combine the experimental method 
with the branch-and-bound method,  there remain many issues to be studied.
In particular, 
the initial setting of the family of 
valid inequalities of the form 
$\sum_{i \in I} x_i \geq \alpha^0_I \ (I \in \IC^0)$ (Section 4.3.1)
should be designed more carefully.
Another important issue is to investigate how to effectively utilize 
the optimal solution information of 
the conic relaxation problem (Section 4.3.2). 
In addition, extensive 
numerical experiment is necessary.

The authors' future interests include applying SCIT to the quadratic assignment problem
(QAP), which is known to be one of the most difficult combinatorial problems. 
They have  participated in the joint project for solving large scale QAPs by the branch-and-bound 
method. See \cite{FUJII2021} for an intermediate report on the project. 
For the lower bounding procedure, the Lagrangian doubly nonnegative (DNN) 
relaxation \cite{KIM2013,ARIMA2017,ARIMA2018} 
and the Newton-bracketing method \cite{KIM2019a,KIM2020}, which were used 
in the numerical results reported in Section 4.5, have been employed in the project. 
For the first time, tai30a and sko42 from QAPLIB \cite{CORAL,ANJOS2021} were solved using their method.  
Although there still remain many unsolved instances in QAPLIB,
 approximate optimal solutions 
which are likely to be optimal are known in all of those instances. The additional cutting 
inequalities discussed in Section 4.3.2 are expected to work effectively to prove that they are 
truly optimal. The branching rules mentioned in Remark~\ref{remark:branching}  can be also used to prove their optimality.  



\end{document}